\documentclass{svmult}
\usepackage{amsfonts,amsmath,amssymb}
\usepackage{graphics, epsfig}
\usepackage{color}
\usepackage{appendix}
\usepackage{verbatim}
\usepackage{ulem}
\usepackage[makeroom]{cancel}
\usepackage[margin=1.5in]{geometry}

 \usepackage[usenames,dvipsnames]{pstricks}
 \usepackage{pst-grad} 
 \usepackage{pst-plot} 
\allowdisplaybreaks

\renewcommand{\thesection}{\arabic{section}}

\numberwithin{equation}{section}
\numberwithin{theorem}{section}
\numberwithin{proposition}{section}
\numberwithin{lemma}{section}
\numberwithin{remark}{section}
\newcommand{\osc}{\operatornamewithlimits{osc}}

\newcommand{\intl}{\int\limits}

\def\Xint#1{\mathchoice
    {\XXint\displaystyle\textstyle{#1}}%
    {\XXint\textstyle\scriptstyle{#1}}%
    {\XXint\scriptstyle\scriptscriptstyle{#1}}%
    {\XXint\scriptscriptstyle\scriptscriptstyle{#1}}%
    \!\int}
\def\XXint#1#2#3{\setbox0=\hbox{$#1{#2#3}{\int}$}
    \vcenter{\hbox{$#2#3$}}\kern-0.5\wd0}
\def\bint{\Xint-}
\def\dashint{\Xint{\raise4pt\hbox to7pt{\hrulefill}}}
\def\dashiint{\bint\kern-0.15cm\bint}

\newcommand{\ovl}[3]{\int_{#1}^{#2}\kern-#3pt\raise4pt\hbox to7pt{\hrulefill}\ }

\newcommand{\ovll}[3]{\intl_{#1}^{#2}\kern-#3pt\raise4pt\hbox to7pt{\hrulefill}\ }

\newcommand{\tvl}[2]{\iint_{#1}\kern-#2pt\raise4pt\hbox to7pt{\hrulefill}\ }

\newcommand{\bye}{\end{document}}

\newcommand{\tvls}[2]{\iint_{#1}\kern-#2pt\raise4pt\hbox to15pt{\hrulefill}\ }

\def\R{\mathbb{R}}
\def\N{\mathbb{N}}
\def\dive{\mathrm{div}}

 \def\B{\mathcal{B}}

\def\L{\mathcal{L}}
\def\C{\mathcal{C}}
\def\T{\mathcal{T}}
\begin{document}
\title*{An Introduction to Barenblatt Solutions for Anisotropic $p$-Laplace Equations.}

\author{Simone Ciani and Vincenzo Vespri}

\institute{Simone Ciani \at Universit\`a degli Studi di Firenze,\\ Dipartimento di Matematica e Informatica "Ulisse Dini", \email{simone.ciani@unifi.it}
\and Vincenzo Vespri \at Universit\`a degli Studi di Firenze,\\ Dipartimento di Matematica e Informatica "Ulisse Dini", 
\email{vincenzo.vespri@unifi.it}}

\maketitle
\abstract{
We introduce Fundamental solutions of Barenblatt type for the equation
\begin{equation} \label{prototype}
    u_t=\sum_{i=1}^N \bigg( |u_{x_i}|^{p_i-2}u_{x_i} \bigg)_{x_i}, \quad \quad p_i >2 \quad \forall i=1,..,N , \quad \quad \text{on} \quad \Sigma_T=\R^N \times[0,T],
\end{equation} \noindent and we prove their importance for the regularity properties of the solutions.}
\vskip.2truecm
\noindent {\bf{MSC 2020}}: 35K67, 35K92, 35B65.
\vskip0.2cm
\noindent
{\bf{Key Words}}: 
Degenerate Orthotropic Parabolic Equations, $p$-Laplace, Anisotropic, Barenblatt Fundamental Solution, Self-Similarity.
\\
{\begin{flushright}
{\it{To celebrate the 60th genethliac of Massimo Cicognani and Michael Reissig.}}
\end{flushright}}

\section*{Introduction.}
Consider the Cauchy problem
\begin{equation} \label{general-anisotropic}
\begin{cases}
u_t =\dive A(x,u,Du), \quad \text{in} \quad \Sigma_T=\R^N \times (0,T),\\
u(x,0)=M \delta(x),
\end{cases}
\end{equation} \noindent where $M>0$, initial datum is the Dirac function $\delta(x)$, the field $A:\Sigma_T \times \R \times \R^N \rightarrow \R^N$ is only measurable and has an anisotropic behavior
\begin{equation} \label{anisotropic-growth}
\begin{cases}
A_i(x,s,z)z_i \ge \Lambda^{*} |z_i|^{p_i}\\
|A_i(x,s,z)| \leq \Lambda_* |z_i|^{p_i-1},
\end{cases}
\end{equation} \noindent
for some constants $\Lambda^*, \Lambda_*>0$ and  $p_i>2$ for any $i \in \{1,..,N\}$. We recall that when all $p_i$s are greater than $2$ the equation is called degenerate. In order to have the existence of solutions, we require the following monotonicity property to the field $A$:
\begin{equation} \label{monotonicity}
    [A(x,s,\xi)-A(x,s,\zeta)] \cdot [\xi-\zeta] >0 , \quad \quad \forall \quad \xi \ne \zeta \quad \text{in} \quad \R^N.
\end{equation} \noindent 
When $p_i \equiv p$ the equation \eqref{general-anisotropic} is named the orthotropic $p$-Laplace, and has nevertheless a different behavior from the classic $p$-Laplace, as its principal part evolves in a way dictated only by the growth in the $i$-th direction. The problem \eqref{general-anisotropic} reflects the modeling of many materials that reveal different diffusion rates along different directions, such as liquid crystals, wood or earth’s crust (see \cite{Ruzicka}). Moreover, as shown in \cite{Mosconi1} the solution to this equation have finite speed of propagation. Note that this is a more reasonable assumption than the usual infinite-speed typical of heat equation, for most of the physical phenomena.
\subsection*{The open problem of regularity.}
The strong nonlinear character and in particular the anisotropy which is prescribed by equation \eqref{general-anisotropic} has proved to be a hard challenge from the regularity point of view. The main difference with standard non linear regularity theory is the growth \eqref{anisotropic-growth} of the operator $A$, usually referred to as {\it{non standard growth}} (see \cite{Antonsev-Shmarev}, \cite{Marcellini-Boccardo}).
This opens the way to a new class of function spaces, called anisotropic Sobolev spaces (see next Section), and whose study is still open and challenging. Even in the elliptic case, the regularity theory for such equations requires a bound on the sparseness of the powers $p_i$. For instance in the general case the weak solution can be unbounded, as proved in \cite{Giaq}, \cite{Marcellini-counter}. However, the boundedness of solutions was proved in \cite{Marcellini-Boccardo} under the assumption that 
\begin{equation} \label{picondition-boundedness}
    \overline{p}<N, \quad \max\{p_1,..,p_N\} < \overline{p}^*,
\end{equation} \noindent where
\begin{equation} \label{def-mediaarmonica e sobolev}
{\overline{p}}:= \bigg(\frac{1}{N} \sum_{i=1}^N \frac{1}{p_i} \bigg)^{-1} , \quad \quad \quad 
{\overline{p}}^*:=\frac{N{\overline{p} }}{N-\overline{p}}.
\end{equation} \noindent 
Regularity properties are proved only on strong assumptions on the regularity of the coefficients (see \cite{Eleuteri-Marcellini-Mascolo},\cite{Marcellini-pq},\cite{Marcellini-nonstandard}). Even in the elliptic case, when the coefficients are rough, H\"older continuity remains still nowadays an open problem.
Indeed, continuity conditioned to boundedness has been proved in \cite{DiBen-Gia-Ves} by means of intrinsic scaling method, but with a condition of stability on the exponents $p_i$ which is only qualitative. Removability of singularities has been considered in \cite{Sk}. We refer to \cite{Eleuteri-Marcellini-Mascolo} and \cite{Mingione}  for a complete survey on the subject and related
bibliography. 
\newpage
\subsection*{Aim of the note.}
We will consider the homogeneous prototype problem
\begin{equation*}
\begin{cases} u_t= \sum_{i=1}^N \bigg( |u_{x_i}|^{p_i-2} u_{x_i} \bigg)_{x_i}, \quad \text{in} \quad \Sigma_T= \mathbb{R}^{N} \times (0,T),\\
u(x,0)= \delta_o.
\end{cases}
\end{equation*} \noindent 
The purpose of this note is to show the importance of a Barenblatt Fundamental solution $\B$ to this equation, paralleling the construction of Fundamental solutions for the $p$-Laplace equation. We will show a fundamental connection between the previous equation and a particular Fokker-Planck equation, as proved for the porous medium equation by \cite{Carrillo-Toscani}. The achievement of such Fundamental solution would provide important tools for the study of regularity of parabolic anisotropic problems as \eqref{general-anisotropic}. As we will see in the sequel, the problem is more delicate than in the isotropic case, because of the lack of radial solutions. In the isotropic case the adoption of radial symmetry brings the equation, set in a proper scale, to a solvable ODE. In the doubly nonlinear case, a non-explicit Barenblatt Fundamental solution has been found with this approach in \cite{Matas-Merker}, using a Leray-Schauder technique. Also in mathematical physics, the use of radial solution is usual. For instance this strategy can be used for the Navier Stokes equation (see \cite{Gurtin}). In our case, as already stated, the anisotropy does not allow the use of radial solutions, and  this fact compels us to look for new ideas.


\subsection*{Acknowledgements.} Both authors wish to acknowledge to Sunra Mosconi his valuable advice. We would like to thank Andrea Dall'Aglio for helpful conversation on the existence for the Cauchy problem with measure data, and Matias Vestberg for information about the construction of Barenblatt solutions for doubly-nonlinear equations. Moreover, both authors are partially founded by INdAM (GNAMPA).

\section*{Preliminaries.}

\subsection*{Self-similar Fundamental solutions. Motivations and historical perspectives.}
The issue of finding Fundamental solutions to elliptic and parabolic equations is one of paramount importance in the study of linear elliptic and parabolic equations (see \cite{DB-PDE}). In nonlinear theory their role is not so evident, and yet the epithet "Fundamental" is iconic, because representation in terms of kernels usually fails. But they are a tool of extraordinary importance in the existence and regularity theory as well as very important to describe the asymptotic behaviour, that's why the name Fundamental Solutions is deserved. Much more information about techniques to be employed, sharp-condition examples and counterexamples can be extracted from the knowledge of a Fundamental solution. A typical example is the Barenblatt Fundamental Solution 
\begin{equation*}
\mathcal{B}(x,t)=  t^{-\frac{N}{\lambda}}\bigg{\{} 1- \gamma_p\bigg( \frac{|x|}{t^{\frac{1}{\lambda}}} \bigg)^{\frac{p}{p-1}} \bigg{\}}_{+}^{\frac{p-1}{p-2}}, \quad t>0,
\end{equation*} for the $p$-Laplace equation 
\begin{equation} \label{plaplaciano}
u_t = \dive(| \nabla u|^{p-2} \nabla u), \quad \text{in} \quad [0,T] \times \R^N, \quad p>1.
\end{equation}
These special solutions can be used to reveal a gap between the elliptic theory and the corresponding parabolic one for $p$-Laplace type equations. Indeed solutions to
\begin{equation} \label{elliptic-plaplace}
    \dive{(|\nabla u|^{p-2} \nabla u)}=0, \quad u \in W^{1,p}_{loc}(\Omega), \quad p>1,
\end{equation} \noindent do obey to a Harnack inequality (see \cite{Serrin}), while the corresponding solutions to the parabolic version of \eqref{elliptic-plaplace} do not in general. We show this briefly. Let $(x_0,t_0)$ be a point of the boundary of the support of $\B$, the free boundary $\{ t=|x|^{\lambda}\}$, and let $\rho>0$. The ball $B_{\rho}(x_0)$ intersects at the time level $t_0-\rho^p$ the support of $x \rightarrow \B(x, t_0-\rho^p)$ in an open set, hence
\[
\B(x_0,t_0)=0, \quad \text{but} \quad \sup_{B_{\rho}(x_0)} \B(x, t_0-\rho^p) >0.
\]
Generalizing the classical heat equation to nonlinear versions, another chief example in evolution theories is the Porous Medium Equation
\begin{equation} \label{PME}
    u_t- \Delta (u^m)=0, \quad m>1.
\end{equation} \noindent This equation, introduced in the last century in connection with a number of physical applications, has been extensively studied (see the monograph \cite{Vazquez-mono}) in parallel to the $p$-Laplace as another prototype of nonlinear diffusive evolution equation, with interest also in the geometry of free boundaries. Fundamental solutions were discovered in $1950$'s by Zeldovich and Kompanyeets in \cite{Zeldo} and Barenblatt \cite{Barenblatt}, and later a complete description has been brought by Pattle in \cite{Pattle}. The discovery of these explicit solutions, usually called Barenblatt solutions since then, has been the starting point of the rigorous mathematical theory  that has been gradually developed since then.\\
The surprising relation between existence and uniqueness of Fundamental solutions and precise asymptotic behaviour relies on the existence of a scaling group under whose action the solutions to the equation are invariant. This implies that a Fundamental solution is self-similar: this is what we call a {\it Barenblatt solution}. Self-similarity has big relevance for the understanding of Fundamental processes in mathematics and physics, as described in \cite{Barenblatt-book}. Self-similar phenomena got in mathematical physics quite early, perhaps with the famous work of Fourier in $1822$ on the analytical theory of heat conduction. In this memoir he performed a construction of a {\it source-type} solution \[ u(x,t)=\frac{A}{\sqrt{t}} f \bigg( \frac{x}{\sqrt{t}} \bigg), \quad \text{for} \quad f(\zeta)= e^{-\frac{\zeta^2}{A}},\,A>0, \] to the heat conduction equation 
\begin{equation} \label{heat}
u_t= \Delta u.
\end{equation} \noindent
Subsequently the phenomena under consideration and their mathematical models became increasingly complicated and very often nonlinear. To obtain self-similar solutions was considered a success in the pre-computer era. Indeed, the
construction of such solutions always reduces the problem to solving the boundary value problems for an ODE, which is a substantial simplification, as we will see in Section \ref{Barenblatt4isotropic}. Furthermore, in 'self-similar' coordinates (as $u \sqrt{t}$, $x/ \sqrt{t}$ for \eqref{heat}), self-similar phenomena become time independent. This enlightens a certain type of stabilization.
Thus during the  pre-computer era, the achievement of a self-similar solution was the only way to understand the qualitative features of the phenomena, and the exponents of the independent variables $x,t$ in self-similar variables were obtained often by dimensional analysis. Dimensional analysis is merely a simple sequence of rules based on the Fundamental covariance principle of physics: all physical laws can be represented in a form which is equally valid for all observers.\\ \\
The very idea of self-similarity is connected with the { \it group of transformations} of solutions (see \cite{Barenblatt-book2}). These groups are already present in the differential equations of the process and are determined by the dimensions of the variables appearing in them: the transformations of the units of time, length, mass, etc. are the simplest examples. This kind of self-similarity is obtained by power laws with exponents that are simple fractions defined in an elementary way from dimensional considerations. Such a course of argument has led to results of immense and permanent importance, as the theory of turbulence and the Reynolds number, of linear and nonlinear heat propagation from a point source, and of a point explosion. Moreover it has enlightened the way toward to a nonlinear theory developed by DiBenedetto (\cite{DB}) with the nowadays well-known method of intrinsic scaling (see also \cite{Urbano}).
\subsubsection*{The group of transformations for the $p$-Laplace equation.}
Let us examine the group of transformations under scaling of the $p$-Laplace equation
\[
u_t= \dive{(|\nabla u|^{p-2}\nabla u)}.
\]
We apply the following dilation in all variables
\[
u'=Ku,\quad \quad x'=Lx, \quad \quad t'=Tt,
\] and impose that the function $u'$ so defined 
\begin{equation} \label{transf}
u'(x',t')=Ku \bigg(\frac{x'}{L}, \frac{t'}{T} \bigg),
\end{equation} is again a solution to the $p$-Laplace equation above. Then by the simple calculations
\[
u_{t'}= \frac{K}{T} u_t \bigg( \frac{x'}{L}, \frac{t'}{T}  \bigg), \quad |\nabla_{x'} u'| = \frac{K}{L} | \nabla u|
\] we arrive to the conclusion that $u'$ is a solution to the $p$-Laplace equation if and only if
\[
TK^{p-2}=L^p.
\] So we obtain a two-parametric transformation group $\mathcal{T}(L,T)$ acting on the set of solutions of the $p$-Laplace equation:
\begin{equation} \label{isotropic-group-transformations}
    (\mathcal{T}u)(x,t)=  \bigg(\frac{L^p}{T} \bigg)^{\frac{1}{p-2}}  u \bigg( \frac{x}{L}, \frac{t}{T}  \bigg).
\end{equation} \noindent and we can conclude the following Lemma.
\begin{lemma}
If $u$ is a solution to the $p$-Laplace equation in a certain class of solutions $\mathcal{S}$ which is closed under dilation in $x,t,u$, then $(\mathcal{T}u)$ given by \eqref{isotropic-group-transformations} is again a solution to the equation in the same class $\mathcal{S}$.
\end{lemma} \noindent
Those special solutions that are themselves invariant under the scaling group are called {\it self similar-solutions}: this means that $(\T u)(x,t)=u(x,t)$ for all $(x,t)$ in the domain of definition, which has to be itself scale-invariant.\\ \\ Suppose now that we have an important information, such as \eqref{Linfty-decay} or conservation of mass. We want to use some of the free parameters to force $\T$ to preserve this important behaviour of the orbit. Analytically it consists in imposing a new relation between two independent parameters, as $K$ and $L$ for instance, and in reducing the transformation to a one-parameter family of scaled functions. Thus we set
\begin{equation} \label{scaling1}
K=L^{-\chi}    ,
\end{equation} \noindent and consequently
\[
K=T^{-\alpha}, \quad \quad \quad L=T^{\beta},
\] with $\alpha, \beta, \chi$ linked by conserving the equation:
\[
\alpha(p, \chi)=\frac{\chi}{\chi(p-2)+2}, \quad \quad \beta(p,\chi)=\frac{1}{\chi(p-2)+2} , \quad \text{unless} \quad \chi=\frac{-2}{(p-2)}.
\] Observing that $\chi= \alpha/\beta$, the equation changes into
\begin{equation}
( \T u)(x,t)=T^{-\alpha} u(x/T^{\beta},t/T)    ,
\end{equation} \noindent where $\alpha, \beta$ are linked by
$\alpha(p-2)+\beta=1$. The condition of preserving the initial mass is
\begin{equation} \label{mass-conservation}
    \int_{\R^N} K u_0\bigg(\frac{x}{L} \bigg) dx = \int_{\R^N} ( \T u_0)(x) dx = \int_{\R^N} u_0(x) dx
\end{equation} \noindent which obliges $KL^N=1$, so that the one parameter family $\T$ will be given by
\begin{equation}
\alpha=\frac{N}{N(p-2)+2}, \quad \quad \beta=\frac{1}{N(p-2)+2} ,\quad p>2.
\end{equation} \noindent 
Observe the formula for the transformation of the initial data (which obviously must satisfy the same transformation) must be
\begin{equation}\label{transf-data}
    (\T u_0)(x)=T^{-\frac{N}{\lambda}} u_0\bigg( \frac{x}{T^{\frac{1}{\lambda}}}  \bigg), \quad \quad \lambda=N(p-2)+p.
\end{equation} \noindent 
In the case of Barenblatt Fundamental solution \eqref{Barenblatt} the couple $(x,t)$ is fixed as a single variable so that
\begin{equation} \label{supertrans}
u(x,t)=t^{-\alpha}u(xt^{-\beta},1)=t^{-\alpha} F(xt^{-\beta}),
\end{equation} \noindent
where $F(\eta)=u(\eta,1)$ is the {\it profile} of the solution.

\begin{remark} A complete theory of existence and uniqueness for the main equation would allow us to obtain self-similar solutions almost for free. Indeed we can consider the solution to the Cauchy problem for scale invariant data, and then use uniqueness to show that this must be self-similar. Let the initial data for instance be of the form
\[
u'(x)=\frac{G(\xi)}{|x|^{\chi}}, \quad \chi \in \R,\, \, \xi= \frac{x}{|x|}, \quad \text{and} \quad  G:\mathbb{S}^{N-1} \rightarrow \R.
\] 
Let us suppose that we are able to solve with uniqueness the Cauchy problem for our equation with this initial data, say the solution is $u$. We produce another solution to the same equation by $\T(u)$ given by \eqref{transf} and if $K=L^{-\chi}$ then the transformed initial data is the same one:
\[
(\T u)(x,0)= K G(\xi) \bigg|\frac{x}{L}\bigg|^{-\chi}=u(x,0)
\] \noindent and so $u$ and $\T(u)$ solve the same Cauchy problem and $u$ is self-similar.
\end{remark}

\subsection*{Notation and settings.}
Given $ {\bf{p}}:= (p_1,..,p_N)$, ${\bf{p}}>1$ with the usual meaning, we assume that the harmonic mean is smaller than the dimension of the space variables
\begin{equation} \label{harmonic-mean}
{\overline{p}}:= \bigg(\frac{1}{N} \sum_{i=1}^N \frac{1}{p_i} \bigg)^{-1}<N,
\end{equation} \noindent and we define the Sobolev exponent of the harmonic mean $\overline{p}$,
\begin{equation} \label{sobolev-exponent} {\overline{p}}^*:=\frac{N{\overline{p} }}{N-\overline{p}}.\end{equation} \noindent We will suppose without loss of generality along this note that the $p_i$s are ordered increasingly. Next we introduce the natural parabolic anisotropic spaces. Given $T>0$ and a bounded open set $\Omega \subset \R$ we define
\[ W^{1,{\bf{p}}}_o(\Omega):= \{ u \in W^{1,1}_o(\Omega) |\,  D_i u \in L^{p_i}(\Omega) \} \]
\[W^{1,{\bf{p}}}_{loc}(\Omega):= \{ u \in L^1_{loc}(\Omega) |\,  D_i u \in L^{p_i}_{loc}(\Omega) \} \]

\[ L^{{\bf{p}}}(0,T;W^{1,{\bf{p}}}_o(\Omega)):= \{u \in L^1(0,T;W^{1,1}_o(\Omega))|\, D_iu \in L^{p_i}(0,T;L^{p_i}_{loc}(\Omega))   \} \]
\[ L^{{\bf{p}}}_{loc}(0,T;W^{1,{\bf{p}}}_o(\Omega)):= \{u \in L^1_{loc}(0,T;W^{1,1}_o(\Omega))|\, D_iu \in L^{p_i}_{loc}(0,T;L^{p_i}_{loc}(\Omega))   \} \]
Now let $A$ be a measurable vector field satisfying the growth conditions \eqref{anisotropic-growth}. By a {\it local weak solution} of \[
u_t=\dive A(x,u,Du), \quad (x,t) \in \Sigma_T,
\]
we understand a function $u \in C^0_{loc}(0,T; L^2_{loc}(\R^N)) \cap L^{\bf{p}}_{loc}(0,T;W^{1,{\bf{p}}}(\R^N))$ such that for all $0<t_1<t_2<T$ and any test function $\varphi \in C^{\infty}_{loc}(0,T;C_o^{\infty}(\R^N))$ satisfies
\begin{equation} \label{anisotropic-localweaksolution}
\int u \varphi \, dx \bigg|_{t_1}^{t_2}+ \int_{t_1}^{t_2} \int (-u \, \varphi_t +A(x,u,Du)\cdot D\varphi) \, dx dt=0,
\end{equation} \noindent where the integral is assumed to be in $\R^N$ when no domain has been specified. By a density and approximation argument this actually holds for any test function of the kind $ \varphi \in W^{1,2}_{loc}(0,T;L^2_{loc}(\R^n))\cap L^{\bf{p}}_{loc}(0,T;W^{1,{\bf{p}}}_o(\Omega))$ for any semirectangular domain $\Omega \subset \subset \R^N$ (see (\cite{Haskovec-Schmeiser}) for a discussion on anisotropic embeddings and semirectangular domains).
\begin{remark} We further give the definition of solution to the prototype equation \eqref{prototype} with $L^1$ initial data, to be used during the development of our work.\\ \\
A measurable function $(x,t) \rightarrow u(x,t)$ defined in $\Sigma_T$ is a { \it weak solution} to the Cauchy Problem \eqref{general-anisotropic} { \it with $L^1$ initial data} if for every bounded open set $\Omega \subset \R$, if
\[ u \in C(0,T; L^1(\Omega)) \cap L^{{\bf{p}}}(0,T; W^{1,{\bf{p}}}(\Omega)), \quad \text{and} \]
\begin{equation} \label{anisotropic-weaksoution}
\begin{aligned}
\int_{\Omega} u(x,t) \varphi(x,t) dx + \int_0^t \int_{\Omega} &\{-u\varphi_t + \sum_{i=1}^N|D_i u|^{p_i-2}D_iuD_i\varphi \} dx d\tau\\
&= \int_{\Omega} u_0(x) \varphi(x,0) dx,
\end{aligned}
\end{equation} \noindent for all $0<t<T$ and all test functions $\varphi \in C^{\infty}(0,T;C^{\infty}_o(\Omega))$.\newline
Weak subsolutions (resp. supersolutions) are defined as above except that in \eqref{anisotropic-weaksoution} equality is replaced by $ \leq $ (resp. $\ge$) and test functions $\varphi \ge 0$ are taken to be nonnegative. 
\end{remark}

\section*{A Self-Similar Solution to the $p$-Laplace Equation.} \label{Barenblatt4isotropic}
Consider the equation
\begin{equation} \label{p-laplace}
\begin{cases}
    u \in C_{loc}(0,T; L^2_{loc}(\R^N)) \cap L^p_{loc} (0,T;W^{1,p}_{loc}(\R^N),\\
    u_t- \dive ( | \nabla u|^{p-2} \nabla u)=0, \quad \text{in} \quad \Sigma_T= \mathbb{R}^N \times (0,T).
    \end{cases}
\end{equation}
In this case we recover the classic $p$-Laplace equation, and we can write explicitly its self-similarity source-solution since the work of Barenblatt \cite{Barenblatt} as
\begin{equation} \label{Barenblatt}
\mathcal{B}(x,t)=  t^{-\frac{N}{\lambda}}\bigg{\{} 1- \gamma_p\bigg( \frac{|x|}{t^{\frac{1}{\lambda}}} \bigg)^{\frac{p}{p-1}} \bigg{\}}_{+}^{\frac{p-1}{p-2}}, \quad t>0
\end{equation} \noindent
with 
\begin{equation} \label{gammap} \lambda=N(p-2)+p,  \quad \quad \gamma_p= \bigg(\frac{1}{\lambda} \bigg)^{\frac{1}{p-1}} \frac{p-2}{p}.
\end{equation} \noindent
We observe that $\B$ satisfies the self-similar transformation \eqref{supertrans}. This function $\mathcal{B}$ solves the Cauchy problem
\begin{equation} \label{CPBarenblatt}
\begin{cases}
u_t- \dive (|\nabla u|^{p-2} \nabla u)=0, \quad \text{in} \quad \mathbb{R}^{N} \times (0, \infty), \\
\mathcal{B}( \cdot, 0)= M \delta_o,
\end{cases}
\end{equation} \noindent where $\delta_o$ is the Dirac measure concentrated at the origin and for every $t>0$ the mass  $M= \| \mathcal{B}(\cdot,t) \|_{L^1(\mathbb{R}^N)}$ is conserved. The initial datum is taken in the sense of measures, which is, for every $\varphi \in C_o(\mathbb{R}^N)$
\[
\int_{\R^N} \B(x,t) \varphi \, dx \rightarrow M \varphi(0), \quad \text{as} \quad t \downarrow 0.
\]
\noindent 
For $t>0$ and every $\rho>0$ we have the important bound
\begin{equation} \label{Linfty-decay} \|\B(\cdot, t) \|_{L^{\infty}(K_{\rho})} = t^{-\frac{N}{\lambda}}, \end{equation} being $K_{\rho}$ the cube of edge $\rho$.
The explicit function $\B$ is classically named Fundamental solution in literature, because it converges pointwise in $\Sigma_T$ to the heat kernel $\Gamma(x,t)$ when $p$ approaches 2,
\[ \B(x,t) \rightarrow (4 \pi)^{N/2}\Gamma(x,t)= \frac{1}{t^{N/2}} e^{-\frac{|x|^2}{4t}}, \quad \text{if} \quad p \downarrow 2, \] but the name does not refer to the kernel property i.e. solutions to \eqref{p-laplace} are not representable as convolutions of $\B$ with initial data. 
Nevertheless all non-negative solutions to \eqref{p-laplace} behave as $t \downarrow 0$ like the Fundamental solution $\B$, and as $|x| \rightarrow \infty $ they grow no faster than $|x|^{p/(p-2)}$. 
Barenblatt Fundamental solutions $\B$ are useful, together with the comparison principle, for proving an intrinsic Harnack estimate (see further Section \ref{Harnack}), uniqueness in existence with $L^1$ data (as in \cite{Kamin-Vazquez}), and more generally to understand the behavior of solutions from the point of view of the physics. In this way, a suitable revisiting of the linear theory had been shaped to face nonlinear equations as the $p$-Laplace. It is possible to build Barenblatt Fundamental solutions centered in $\bar{x}$ with initial datum at a time $\bar{t}$ in the following way
\begin{equation} \label{generalbarenblatt}
\B_{k,\rho} (x,t,\bar{x},\bar{t})= \frac{k \rho^N}{S^{\frac{N}{\lambda}}(t)} \bigg{\{} 1- \bigg( \frac{|x-\bar{x}|}{S^{\frac{1}{\lambda}}(t)} \bigg)^{\frac{p}{p-1}} \bigg{\}}_{+}^{\frac{p-1}{p-2}}, \quad \quad \quad\lambda=N(p-2)+p,
\end{equation} \noindent
with 
\begin{equation} \label{S}
S(t)= \lambda \bigg( \frac{p}{p-2} \bigg)^{p-1} k^{p-2} \rho^{N(p-2)} (t-\bar{t}) + \rho^{\lambda}.
\end{equation} \noindent
These functions enjoy the following important properties.
 \begin{enumerate}
  \item They are weak solutions to \eqref{p-laplace} in $\mathbb{R}^{N} \times \{ t> \bar{t} \}$. 
  
  \item If we fix $t=\bar{t}$ then $B_{k,\rho} \equiv 0 $ for all $x \in \bigg( \mathbb{R}^N - B_{\rho}(\bar{x}) \bigg)$ and for $t>\bar{t}$ the function $x \rightarrow B_{k,\rho}$ vanishes, in a $C^1$ fashion, across the boundary of the ball $\{|x-\bar{x}| < S^{\frac{1}{\lambda}}(t) \}$.\\ Their support evolves compactly:
  \begin{equation} \label{supportB}
  \text{supp} \bigg( B_{k,\rho} (x,t,\bar{x},\bar{t}) \bigg) = \bigg{\{} |x-\bar{x}| \leq S^{\frac{1}{\lambda}}(t) \bigg{\}} \times [\bar{t}, t^{*}],
  \end{equation} \noindent thus
  \begin{equation}
  \text{supp} \bigg( B_{k,\rho} (x,t,\bar{x},\bar{t}) \bigg) \subseteq B_{S^{1/\lambda}(t^{*})} (\bar{x}) \times [\bar{t}, t^{*}].
  \end{equation}

\item They are bounded for fixed $\rho$ and $k \in \mathbb{R}^{+}$:
\begin{equation} \label{boundednessB}
    B_{k,\rho} (x,t,\bar{x},\bar{t}) \leq k, \quad \quad \quad x \in \mathbb{R}^N.
\end{equation}
\end{enumerate} \noindent 
In the sequel when no explicit formula for a solution as \eqref{generalbarenblatt} (as in \eqref{prototype}), we will refer to a Barenblatt Fundamental Solution as a function (resp. to \eqref{prototype}) satisfying properties analogous to 1-3 above.

\subsubsection*{The construction of $\B$: reduction to an isotropic Fokker-Planck equation.}
As far as we know if we look for a Barenblatt Fundamental solution as $\B$, we have to impose the condition \eqref{Linfty-decay}, because this is the behaviour that non-negative solutions to the $p$-Laplace Cauchy problem with the right decay of the initial datum do satisfy (see \cite{DB} Theorem 4.5). This motivates us to apply the following (formal) transformations to the equation \eqref{p-laplace}
and 
\begin{equation} \label{iso-transform}
    \begin{cases}
    u(x,t)=t^{-\frac{N}{\lambda}} v(xt^{\alpha},t)=v(y,t),\\
    y=x t^{\alpha}, \quad \quad \alpha=-\frac{1}{\lambda},
    \end{cases} \Rightarrow 
    \begin{cases}
    u_x= t^{-\frac{N}{\lambda}} v_y y_x= t^{\alpha-\frac{N}{\lambda}} v_y,\\
    \frac{\partial}{\partial x}= t^{\alpha} \frac{\partial}{\partial y}.
    \end{cases}
\end{equation} \noindent
\begin{remark} \label{backtosimilarity}
We notice that the applied transformation does not belong to the group of transformations \eqref{isotropic-group-transformations}, so we expect that equation \eqref{p-laplace} turns into another one. This is what is called in \cite{Vazquez-whole-space} the continuous rescaling: as the change of variables \eqref{iso-transform} belongs to the transformation group only for the fixed time $t=1$, source-type solutions transform into stationary profiles of the transformed equation. \end{remark} \noindent 
By direct calculation we obtain
\[ \begin{aligned}
        u_t & = -\frac{N}{\lambda} t^{-\frac{N}{\lambda}-1}v+ t^{-\frac{N}{\lambda}} \bigg[\sum_{i=1}^{N} v_{y_i} (y_i)_{t} +  v_t \bigg]=\\
        &-\frac{N}{\lambda} t^{-\frac{N}{\lambda}-1}v+ t^{-\frac{N}{\lambda}} \bigg[ \nabla_y v \cdot \frac{\alpha y}{t}+  v_t \bigg]
    \end{aligned} \]\noindent and 
 \begin{equation}
        \nabla_x u= t^{\alpha-\frac{N}{\lambda}} \nabla_y v.
    \end{equation} \noindent We set
    \begin{equation}
       \tilde{v} (y,\tilde{t})= \tilde{v}(y,\ln (t))= v(y,t), \quad \Rightarrow \quad \tilde{v}_t= \tilde{v}_{\tilde{t}} t^{-1}= v_t
    \end{equation} \noindent and the equation \eqref{p-laplace} becomes, by multiplying it for $t^{\frac{N}{\lambda}+1}$
\[       
\begin{aligned} 
 \tilde{v}_{\tilde{t}}= & \frac{N}{\lambda} v- \frac{N}{\lambda} \nabla_y \tilde{v} \cdot y + t^{\alpha} \nabla_y \cdot \bigg[ t^{(\alpha-\frac{N}{\lambda})(p-1)} |\nabla_y \tilde{v} |^{p-2} \nabla_y \tilde{v} \bigg] t^{\frac{N}{\lambda}+1}=\\
&\frac{N}{\lambda} v- \frac{N}{\lambda} \nabla_y \tilde{v} \cdot y + \nabla_y \cdot \bigg[|\nabla_y \tilde{v} |^{p-2} \nabla_y \tilde{v} \bigg] t^{\alpha+(\alpha-\frac{N}{\lambda})(p-1)+\frac{N}{\lambda}+1}=\\
& \frac{N}{\lambda} v- \frac{N}{\lambda} \nabla_y \tilde{v} \cdot y +\nabla_y \cdot \bigg[|\nabla_y \tilde{v} |^{p-2} \nabla_y \tilde{v} \bigg],
\end{aligned}
    \] \noindent being $\alpha=-\frac{1}{\lambda}$. So we obtain the isotropic Fokker-Planck equation
\begin{equation} \label{Fk-Pk-isotropa}
\tilde{v}_{\tilde{t}}= \nabla_y \cdot \bigg( | \nabla_y \tilde{v}|^{p-2} \nabla_y \tilde{v} + \frac{y \tilde{v}}{\lambda}   \bigg).
\end{equation} \noindent

\subsubsection*{Barenblatt solution solves the isotropic Fokker Planck equation.}
Consider the Barenblatt function $\B(x,t)$, with explicitly scaled space variables
\begin{equation} 
\B(x,t)= t^{-\frac{N}{\lambda}} \bigg{\{} 1-\gamma_p \bigg( {\sqrt{\sum_{i=1}^{N} \bigg(\frac{x_i}{t^{\frac{1}{\lambda}}} \bigg)^2}}  \bigg)^{\frac{p}{p-1}}    \bigg{\}}_{+}^{\frac{p-1}{p-2}}.
\end{equation} \noindent
We claim that $\B$ solves the stationary version of \eqref{Fk-Pk-isotropa}, by taking the flux to be zero, i.e.
\[
 | \nabla_y \tilde{v}|^{p-2} \nabla_y \tilde{v} + \frac{y \tilde{v}}{\lambda} =0.
\]
We have, by setting $y_i= x_i t^{-\frac{1}{\lambda}}$, that 
\[ \B(y,t)=t^{-\frac{N}{\lambda}} \bigg{\{} 1-\gamma_p |y|^{\frac{p}{p-1}}    \bigg{\}}_{+}^{\frac{p-1}{p-2}} = t^{-\frac{N}{\lambda}} \bigg{\{} 1-\gamma_p \bigg( \frac{\sqrt{\sum_{i=1}^{N} x_i^2}}{t^{\frac{1}{\lambda}}}  \bigg)^{\frac{p}{p-1}}    \bigg{\}}_{+}^{\frac{p-1}{p-2}}= \B(x,t)
\]
and thus the function
\[ \mathcal{C}(y,t)=\bigg{\{} 1-\gamma_p |y|^{\frac{p}{p-1}}    \bigg{\}}_{+}^{\frac{p-1}{p-2}}\] \noindent is independent from $t$, and 
\[ \begin{aligned}
\nabla_y \C  = & - \gamma_p \bigg( \frac{p}{p-2} \bigg) \bigg{\{} 1-\gamma_p |y|^{\frac{p}{p-1}}    \bigg{\}}_{+}^{\frac{1}{p-2}} |y|^{\frac{2-p}{p-1}} y =\\
& = - \gamma_p  \bigg( \frac{p}{p-2} \bigg)\C ^{\frac{1}{p-1}} |y|^{\frac{2-p}{p-1}} y.
\end{aligned}
\]
Thus by calculation we have that $\C(y)= t^{\frac{N}{\lambda}}\B(y,t)$ solves the zero flux equation
\[
\begin{aligned}
| \nabla_y \C |^{p-2} \nabla_y \C  &+ \frac{y \C }{\lambda}=\\
&\bigg[ \gamma_p \bigg( \frac{p}{p-2} \bigg) \C ^{\frac{1}{p-1}}   \bigg]^{p-2} |y|^{\frac{2-p}{p-1}(p-2)} |y|^{p-2} \bigg[ -\gamma_p \bigg(\frac{p}{p-1}  \bigg) \C^{\frac{1}{p-1}} |y|^{\frac{2-p}{p-1}}y \bigg]+ \frac{y \C}{\lambda}=\\
& \C \bigg[\frac{1}{\lambda}-\gamma_p \bigg( \frac{p}{p-2} \bigg)^{p-1}    \bigg] y=0, \quad \text{for} \quad  \gamma= \bigg( \frac{p-2}{p} \bigg)^{p-1}\frac{1}{\lambda}.
\end{aligned}
\] 
Consequently, so does $\B(x,t)$. Now we show that the converse reasoning holds too, in order to show how the whole calculation is in fact reduced to a ODE solution.
\subsubsection*{Function $\C$ solves a particular ODE.}
Consider
\begin{equation} \label{1-varBarenblatt}
\C (\eta)=\bigg{\{} 1-\gamma_p \eta ^{\frac{p}{p-1}}    \bigg{\}}_{+}^{\frac{p-1}{p-2}}= \C (|y|), \quad \eta >0.
\end{equation} \noindent
In $0 \leq \eta < \bigg( \frac{1}{\gamma_p} \bigg)^{\frac{p-1}{p}}$ we have
\[ \C (\eta)^{\frac{p-2}{p-1}} =1-\gamma_p \eta^{\frac{p}{p-1}}. \]
We derive the equation to obtain
\[
\bigg( \frac{p-1}{p-2} \bigg) \C (\eta)^{-\frac{1}{p-1}} \C'(\eta) d\eta = - \bigg( \frac{p-2}{p} \bigg) \frac{1}{\lambda^{1/(p-1)}} \bigg( \frac{p}{p-1}  \bigg) \eta^{\frac{1}{p-1}} d\eta.
\]
Now, we manipulate the equation with $\C'(\eta) \leq 0$, because
\[
\C '(\eta)= \bigg( \frac{p-1}{p-2} \bigg) \bigg{\{}1-\gamma_p \eta^{\frac{p}{p-1}}   \bigg{\}}_+^{-\frac{1}{p-2}} \bigg( -\gamma \bigg( \frac{p}{p-1} \bigg) \eta^{\frac{1}{p-1}} \bigg) \leq 0
\] 
so that
\[
\bigg( \frac{(-\C'(\eta))^{p-1}}{\C(\eta)}  \bigg)^{\frac{1}{p-1}}= \bigg( \frac{\eta}{\lambda} \bigg)^{\frac{1}{p-1}}
\]
and so the desired mono-dimensional Fokker-Planck equation is obtained
\begin{equation} \label{monoFokker}
|\C'(\eta)|^{p-2} \C'(\eta) + \frac{\eta \C(\eta)}{\lambda}=0.
\end{equation} \noindent 
If one reads conversely from the end to the beginning of these calculations, it is clear how to arrive to a solution to the isotropic Fokker Planck equation \eqref{Fk-Pk-isotropa} by imposing radial symmetry.


\section*{Solving the isotropic Cauchy Problem with measure data.} \label{SolvingCP}
Suppose now that we are not able to solve by radial symmetry the isotropic Fokker-Planck equation \eqref{Fk-Pk-isotropa}. If we look for a solution to \eqref{CPBarenblatt} that exhibits the properties \eqref{supportB}-\eqref{boundednessB}, we may adopt the following strategy. First we find a general solution $u$ to \eqref{CPBarenblatt} with datum the Dirac measure $\delta_o$, we show that it is positive by the maximum principle, and then we use the transformation \eqref{iso-transform} to get a solution $w$ to \eqref{Fk-Pk-isotropa}. Observe that a comparison principle for subsolutions to the $p$-Laplace equation can be transported to a comparion principle for subsolutions to the isotropic Fokker-Planck equation. But we need a solution to the stationary Fokker-Planck equation to recover the self-similarity (see Remark \ref{backtosimilarity}), so that we can control the behavior for all times by scaling, and we gain for free the correct evolution of its support.
More generally speaking, if the initial data in \eqref{anisotropic-weaksoution} is given by
\begin{equation}
u_0(\cdot,0)= \mu,    
\end{equation} \noindent where $\mu$ is a $\sigma$- finite Borel measure in $\R^N$, then we say that $u$ is a {\it weak solution} of \eqref{anisotropic-weaksoution} {\it with measura data} if for every bounded open set $\Omega \subset \R^N$ and $ \forall t \in (0,T)$, $u$ satisfies the above integral equality \eqref{anisotropic-weaksoution} with the right-hand side replaced by
\[ \int_{\Omega} \varphi(x,0) d \mu,   \]
$\forall \varphi \in C^1( \overline{\Omega_T})$ such that $ x \rightarrow \varphi(x,t)$ is compactly supported in $\Omega \, \forall t \in [0,T]$. \\ \\
In the pioneering work \cite{DB-He} for the isotropic $p$-Laplace, the authors consider a way of measuring the growth of a function $f \in L^1_{loc}(\R^N)$ as $|x| \rightarrow \infty$ by setting
\[
|\| f \||_r:= \sup_{\rho \ge r} \rho^{-\lambda /(p-2)} \int_{B_{\rho}} |f| dx, \quad r>0, \quad  \lambda=N(p-2)+p.
\]
Note that if $f \in L^1(\R^N)$ then $|\| f \||_r < \infty,\,  \forall r>0 $. Similarly, if $\mu$ is a $\sigma$-finite Borel measure in $\R^N$, we set
\[
|\| \mu \||_r:= \sup_{\rho \ge r} \rho^{-\lambda /(p-2)} \int_{B_{\rho}} |d \mu|,
\]
where $|d \mu|$ is the variation of $\mu$.\\
In that Fundamental work, the authors demonstrate the existence of a weak solution to the problem \eqref{anisotropic-weaksoution} in its isotropic configuration, within $\Sigma_T=\Sigma_T(\mu)$, where 
\begin{equation} \label{isotropic-time}
T(\mu)= \begin{cases}  
C_0(N,p) \bigg[ \lim_{r \rightarrow \infty} |\| \mu \||_r  \bigg]^{(2-p)}, & \text{if} \quad \lim_{r \rightarrow \infty} |\| \mu \||_r >0\\
+ \infty &\text{if} \quad \lim_{r \rightarrow \infty} |\| \mu \||_r =0.
\end{cases}
\end{equation} \noindent So the existence is proved in a cylindrical domain whose last time $T$ is dictated by the behavior at infinity of the initial measure $\mu$. The method relies on suitable estimates and compactness, which permit a standard limiting process. Indeed, given a $\sigma$- finite Borel measure $\mu $ in $\R^N$ satisfying $|\|\mu \||_r < \infty $ for some $r>0$, there exists a sequence of regular functions $\{ u_{0,n} \}_{n \in \N} \in C_o^{\infty}(\R^N)$ such that $\forall \varphi \in C_o(\R^N)$ we have
\[
\int_{\R^N} u_{0,n} \varphi \, dx \rightarrow \int_{\R^N} \varphi d \mu, \quad {\&}\quad |\| u_{0,n} \||_r \rightarrow  |\| \mu\||_r, \quad r>0. 
\]
The Cauchy Problem 
\begin{equation} \label{CP-isotropic+bounded-domain} \begin{cases} 
u_t- \dive (|Du|^{p-2}Du)=0 & \text{in} \quad \Sigma_T, \quad p>2, \\
u(\cdot,0)= u_{0,n}.
\end{cases}
\end{equation} \noindent
has a unique solution $u_n$, global in time (see \cite{Benilan-Crandall}). Next, the authors prove the following estimates, for all $0<t< T_r(\mu):= C_0 [\|| \mu \||_r]^{(2-p)}$, $\forall \rho \ge r >0$:
\begin{equation} \label{isoestimate1}
\| |u (\cdot, t) \||_r \leq C_1(N,p) \| | \mu \||_r,
\end{equation}

\begin{equation} \label{isoestimate2}
\| u(\cdot,t) ||_{L^{\infty}(B_{\rho})} \leq C_2(N,p) t^{-N/ \lambda} \rho^{p/(p-2)} |\| \mu \||_r^{p/\lambda},
\end{equation}

\begin{equation} \label{isoestimate3}
\|Du(\cdot,t)||_{L^{\infty}(B_{\rho})} \leq C_3(N,p) t^{-(N+1)/\lambda} \rho^{2/(p-2)} \|| \mu \||_r^{2/\lambda},
\end{equation}
\begin{equation}
    \int_0^t \int_{\Omega} |Du|^q dx d \tau \leq C_4(N,P,\epsilon, \text{diam} \Omega) \,  |\| \mu \||_r^{C_5(N,p,\epsilon)}, \quad q=p-(N+\epsilon)/(N+1),
\end{equation} \noindent and in particular with $\epsilon=1$ we obtain
\begin{equation}
\int_0^t \int_{B_{\rho}} |Du|^{p-1} dx d\tau \leq C_5(N.p) t^{1/\lambda} \rho^{1+\lambda/(p-2)} \, | \| \mu \| |_r^{1+(p-2)/\lambda}
\end{equation}
Moreover the function $(x,t) \rightarrow Du(x,t)$ is H\"older continuous in $\overline{\Omega} \times [ \eta, T(\mu)-\eta], \, 0<\eta<T(\mu)$, with H\"older constants and exponents depending upon $N,p,C_1,..,C_4$, $\text{diam}\Omega, \eta, |\| \mu \||_r$. It can be shown that their estimates are sharp, by means of Barenblatt solutions. Finally, the estimates above \eqref{isoestimate1}-\eqref{isoestimate3} with a monotonicity property as \eqref{monotonicity}, permit to pass to the limit in the approximating problems \eqref{CP-isotropic+bounded-domain}. 


\section*{An application of $\B$ to intrinsic Harnack estimates.} \label{Harnack}
In this section we outline the importance of the construction of a Barenblatt Fundamental solution for the aim of proving regularity. Indeed the rough idea is that once that we have a solution of \eqref{p-laplace} whose support and positivity can be easily manipulated, by means of a comparison argument is possible to expand the positivity set of a whatever solution that is bigger than the Fundamental one in the parabolic boundary. More precisely we will review the proof of the following Theorem of \cite{DB}.
\begin{theorem} \label{harnackpotente}
Let $u$ be a non-negative weak solution of equation \eqref{p-laplace} in $\Omega_T=\Omega \times [0,T]$ where $\Omega \subset \R^N$ bounded open set. Fix a point $(x_0,t_0) \in \Omega_T$ and assume $u(x_0,t_0)>0$. There exist constants $\gamma >1$ and $C>1$
, depending only on $N,p$, such that 
\begin{equation} \label{Harnack-estimate}
u(x_0,t_0) \leq \gamma \inf_{B_{\rho}(x_0)} u( \cdot, t_0+ \theta), \quad \theta= \frac{C \rho^p}{[u(x_0,t_0)]^{p-2}},    
\end{equation} \noindent  provided the cylinder 
\begin{equation}
Q_{4 \rho} (\theta)= \{ |x-x_0| <4 \rho \} \times \{ t_0-4 \theta, t_0+ 4 \theta \}
\end{equation} \noindent is contained in $ \Omega_T$.
\end{theorem}

\begin{remark}
As we can see, the geometry is intrinsically defined by the value of the solution in $(x_0,t_0)$. This brings to light a difficulty in exposition, as a priori weak solutions to \eqref{p-laplace} are not meant to be well defined in every point. Nonetheless by standard regularity theory we know that local weak solutions to \eqref{p-laplace} are locally H\"older continuous, and so they are well defined pointwise as elements of $C(0,T;W^{1,p}_{loc}(\Omega))$. \end{remark}

\begin{remark}
The constants $\gamma$ and C in previous Theorem tend to infinity as $p$ tend to infinity, but they are stable as $p \downarrow  2$ in the following meaning
\begin{equation}
    \lim_{p \downarrow 2} \gamma(N,p)= \gamma(N,p), \quad \text{and} \quad \lim_{p \downarrow 2} C(N,p)= C(N,p).
\end{equation} \noindent 
\end{remark}

\subsection*{Outline of the proof of Theorem \ref{harnackpotente}.}

For the sake of conciseness ad to the aim of highlighting the importance of Barenblatt Fundamental solutions, we will demonstrate only the case when $p$ is not too close to $2$. The proof for $p \in (2,5/2]$ uses local comparison functions built especially to do the same job of $\B$, being subsolutions of \eqref{p-laplace} and observing the same ordering imposed by the following Lemma.
\begin{lemma}
Let $u,v$ be two solutions of \eqref{p-laplace} in $\Omega_T=\Omega \times [0,T]$ such that $u,v \in C(0,T; L^2(\Omega)) \cap L^p(0,T; W^{1,p}(\Omega)) \cap C(\overline{\Omega_T})$. If $u \ge v$ in the parabolic boundary of $\Omega_T$, then $u \ge v$ in $\Omega_T$.
\end{lemma} \noindent 
\noindent {\sc Step 1.} \textit{Transforming the problem by scaling.} \newline
Let $(x_0,t_0) \in \Omega_T$, $\rho>0$ to be fixed a posteriori, assume that $u(x_0,t_0)>0$ and for a constant $C$ to be determined later let $Q_{4\rho}$ be the box
\begin{equation} \label{thebox}
    Q_{4\rho}= \{ |x-x_0| <4 \rho \} \times \bigg{\{} t_0- \frac{4C \rho^p}{[u(x_0,t_0)]^{p-2}}, t_0+ \frac{4C \rho^p}{[u(x_0,t_0)]^{p-2}} \bigg{\}}.
\end{equation} \noindent
Now introduce the change of variables
\begin{equation} \label{change4harnack}
\Phi (x,t)= \bigg(\frac{x-x_0}{\rho}, \,\frac{(t-t_0)[u(x_o,t_0)]^{p-2}}{\rho^p} \bigg), \quad \quad \Phi(Q_{4\rho})= B_4 \times (-4C,4C)=:Q
\end{equation} \noindent 
Let us denote again with $x,t$ the new variables $\Phi(x,t)$, and observe that the function
\begin{equation}
v(x,t)= \frac{1}{u(x_0,t_0)} u \bigg( x_0+ \rho x, \frac{t \rho^p}{[u(x_0,t_0]^{p-2}} \bigg),    
\end{equation} \noindent is a bounded non-negative solution to the Cauchy problem
\begin{equation} \label{CP-rescaled4harnack}
    \begin{cases}
    v_t-\dive(|Dv|^{p-2}Dv)=0,& (x,t) \in Q\\
    v(0,0)=1.
    \end{cases}
\end{equation} \noindent 
Theorem \ref{harnackpotente} will be proved, as shown by a simple converse rescaling, if we are able to find constants $\gamma_o \in (0,1]$, $C>1$ depending only upon $N,p$ holding the inequaity
\begin{equation} \label{harnackpotente-rescaled}
    \inf_{B_1} v(x,C) \ge \gamma_o.
\end{equation} \noindent 
The constant $\gamma_o$ defined successively in \eqref{gammazero} tends to zero as $p \downarrow 2$. \\ \\
\noindent {\sc Step 2.} \textit{Finding qualitatively a point where $v$ equals a power-like function of time.} \newline
We consider the family of nested and expanding boxes
\begin{equation} \label{nested-expanding-boxes}
    Q_{\tau}= \{ |x| < \tau\} \times (-\tau^p,0], \quad \quad \tau \in (0,1]
\end{equation} \noindent and for each of these boxes we consider the numbers
\begin{equation} \label{numbers}
    M_{\tau}= \sup_{Q_{\tau}} v, \quad \quad \quad N_{\tau}= (1-\tau)^{-b},
\end{equation} \noindent where the number $b>0$ will be suitably defined later to render quantitative the following estimate. As $M_0=N_0$ and considering that $M_{\tau}$ remains a bounded function of $\tau$ (because $v$ is a {\it bounded} solution) while $N_{\tau} \rightarrow + \infty $ when $\tau$ tends to $1$, we can choose a number $\tau_o$ to be the largest root of the equation 
\[ M_{\tau} = N_{\tau}.
\] This implies by construction
\begin{equation}
    \sup_{Q_{\tau}} v \leq N_{\tau}, \quad \quad \forall \tau > \tau_o.
\end{equation}\noindent Since $v$ is continuous in $Q$ there exists at least one point $(\bar{x},\bar{t}) \in \overline{Q_{\tau_o}}$ such that
\begin{equation}
    v(\bar{x},\bar{t})= N_{\tau_o}= (1-\tau_o)^{-b}.
\end{equation}

\noindent {\sc Step 3.} \textit{Ordering $v$ and $(1-\tau_o)^{-b}$ within a small ball centered in $\bar{x}$.}\newline
Let
\[
R= \frac{1-\tau_o}{2},
\]
and consider the cylinder $ [(\bar{x},\bar{t})+Q(R^p,R)]= {\{}|x-\bar{x}| < R{\}} \times \{ \bar{t}- R^p, \bar{t} \}.$ As $\tau_o \in (0,1]$ we have the inclusion $[(\bar{x},\bar{t})+Q(R^p,R)] \subset Q_{\frac{1+\tau_o}{2}}$ which implies
\[
\sup_{[(\bar{x},\bar{t})+Q(R^p,R)]} v \leq N_{\frac{1+\tau_o}{2}} = 2^b(1-\tau_o)^{-b}=: \omega.
\] Now we use H\"older continuity of the function $v$ in the fashion of Proposition 3.1 of Chap.III of \cite{DB}, choosing $b>0$ so large that the starting one of the family of shrinking cylinders is contained in $ [(\bar{x},\bar{t})+Q(R^p,R)]$. Hence 
there exist $\gamma >1$ and $a, \varepsilon_o \in (0,1)$ such that for all $r \in (0,R]$ we have
\begin{equation}
   \begin{aligned}
\osc_{[(\bar{x},\bar{t})+Q(R^p,R)]}v(\cdot, \bar{t}) &\leq \gamma ( \omega + R^{\varepsilon_o}) \bigg( \frac{r}{R} \bigg)^{a}\\
& \leq 2^{b+1} \gamma (1-\tau_o)^{-b} \bigg( \frac{r}{R}\bigg)^{a}
\end{aligned}
\end{equation} \noindent
We let $r=\sigma R$ and we choose $\sigma$ so small that for all $\{ |x-\bar{x}|<\sigma R \}$ we obtain
\begin{equation}
    \begin{aligned}
    v(x,\bar{t}) &\ge v(\bar{x},\bar{t})-2^{b+1} \gamma (1-\tau_o)^{-b} \sigma^a\\
    &(1-2^{b+1})\gamma \sigma^a)(1-\tau_o)^{-b}\\
    &\frac{1}{2} (1-\tau_o)^{-b}, \quad \forall \{|x-\bar{x}| \}< \sigma R, \quad R= \frac{1}{2}(1-\tau_o)
    \end{aligned}
\end{equation} \noindent 
\noindent {\sc Step 5.} \textit{Expansion of the positivity set and conclusion.}\newline
In this last step we will choose the constants $b>1$ and $C>1$ so that the qualitative largeness of $v(\cdot,\bar{t})$ in the small ball $B_{\sigma R}(\bar{x})$ turns into a quantitative bound below over the full sphere $B_1$ at some later time level $C$. This will be carried on by means of the comparison with the functions $\B_{k,\rho}$ defined in \eqref{generalbarenblatt} by
\[
\B_{k,\rho} (x,t,\bar{x},\bar{t})= \frac{k \rho^N}{S^{\frac{N}{\lambda}}(t)} \bigg{\{} 1- \bigg( \frac{|x-\bar{x}|}{S^{\frac{1}{\lambda}}(t)} \bigg)^{\frac{p}{p-1}} \bigg{\}}_{+}^{\frac{p-1}{p-2}},
\]
\[
S(t)= \lambda \bigg( \frac{p}{p-2} \bigg)^{p-1} k^{p-2} \rho^{N(p-2)} (t-\bar{t}) + \rho^{\lambda}.
\] Indeed, we choose appropriately
\[
k=\frac{1}{2} (1-\tau_o)^{-b}, \quad \quad \rho= \sigma R,
\] and we observe that at the time level $t=C$ the support of $\B_{k,\rho}(\cdot,C,\bar{x},\bar{t})$ is the ball
\[
|x-\bar{x}|^{\lambda} < S(t)= \{d\gamma^{p-2} (1-\tau_o)^{(N-b)/(p-2)} (C-\bar{t})+ (\sigma R)^{\lambda})   \}
\]for 
\[
\gamma(N,b)=\frac{1}{2} \bigg( \frac{\sigma}{2}  \bigg)^{N}, \quad \text{and} \quad d= \lambda \bigg( \frac{p}{p-2}  \bigg)^{p-1}.
\] Now choose
\begin{equation} \label{gammazero}
b=N, \quad \quad \quad C=\frac{3^{\lambda}}{d \gamma^{p-2}},
\end{equation} so that the support of $\B_{k,\rho}(\cdot, C, \bar{x},\bar{t})$ contains $B_2$ and we can use the comparison principle with $v$ as we have in $B_{\rho}$
\begin{equation}
    v(\cdot,\bar{t}) \ge \frac{1}{2}(1-\tau_o)^{-N}=k\ge B_{k,\rho}(\cdot,\bar{t}).
\end{equation} \noindent Thence 
\begin{equation}
    \begin{aligned}
    \inf_{x \in B_1} v(x,C) &\ge \inf_{x \in B_1} B_{k,\rho}(x, C, \bar{x},\bar{t})\\
    & \ge 2^{-(1+2N/\lambda)} \bigg( \frac{\sigma}{2} \bigg)^{N} \bigg{\{}1 -\bigg( \frac{2}{3} \bigg)^{\frac{p}{p-1}} \bigg{\}}^{\frac{p-1}{p-2}}=: \gamma_o,
 \end{aligned}
    \end{equation} \noindent and the proof is concluded.

\section*{Looking for a Barenblatt-type solution to \eqref{prototype}.}
In this section we calculate the right exponents for the transformation of the equation \eqref{prototype} into an anisotropic Fokker-Planck equation. Next we observe that the impossibility of using radial solutions does not allow us to obtain an ODE from the Fokker-Planck equation. Finally we show a strategy to find a non-explicit Barenblatt Fundamental solution.

\begin{remark} Observe initially that we can construct a source-type solution, but that unfortunately has not a compact support. Indeed, consider the following solution to \eqref{prototype}. Let $i\in \{ 1,..,N \}$ and
\begin{equation} \label{monodimensionalsolutions}
    f_i(x_i,t,T_i)= \kappa_i \bigg( \frac{|x_i|^{p_i}}{(T_i-t)} \bigg)^{\frac{1}{p_i-2}}, \quad \kappa_i=\kappa_i(p_i)>0, \quad p_i>2,
\end{equation} \noindent be solutions of the equations
\begin{equation} \label{monoequation}
u_t-(| u_{x_i}|^{p_i-2} u_{x_i} )_{x_i}=0, \quad x_i \in \R, \quad t>0.
\end{equation} \noindent Then the function
\begin{equation} \label{anisosolution-noncompact}
\begin{aligned}
\mathcal{F}(x,t)&= \sum_{i=1}^N  f_i(x_i,t,T_i)\\
&= \sum_{i=1}^N \kappa_i \bigg( \frac{|x_i|^{p_i}}{(T_i-t)} \bigg)^{\frac{1}{p_i-2}}
\end{aligned} \end{equation} \noindent
solves the prototype equation \eqref{prototype}. The same can be done by choosing  $f_i \equiv \B_i$ the mono-dimensional Barenblatt solutions solving \eqref{monoequation}. These functions reveal some of the physical aspects of equation \eqref{prototype}: for instance they can be used to show that the lifetime of solutions is dictated by the largest exponent $p_N$ in the case of large initial mass (see Remark 3 in \cite{Tedeev-estimates}). Unfortunately solutions so-built do not have a compactly supported evolution and we cannot use them to expand the positivity by comparison as done in step 5 of Section 4.1.
 \end{remark}

\subsection*{Finite speed of propagation.}
Consider the Cauchy problem
\begin{equation} \label{anisocauchy}
\begin{cases}
u_t= \dive ( {\bf{A}}(t,x,u,\nabla u)), & \text{in} \quad \Sigma_T= \mathbb{R}^N \times (0,T),\\
u(x,0)=u_0(x) \in L^2(\mathbb{R}^N),
\end{cases}
\end{equation} \noindent where ${\bf{A}}(t,x,u,\nabla u)= ( A_i(x,t,u,\nabla u))_{i=1,..,N}$ is a Caratheodory vector field satisfying the growth conditions \eqref{anisotropic-growth}.
In \cite{Mosconi1} the authors proved the following decay properties, that will be useful to us to intercept the right exponents in the scaling transformation leading to the Fokker-Planck equation for solutions to \eqref{prototype}.

\begin{theorem} \label{Sunra}
Suppose that $p_i>2$ for all $i \in \{1,..,N\}$.
Let $u$ be a local weak solution to \eqref{anisocauchy} in $\Sigma_T$ under the growth conditions \eqref{anisotropic-growth} with
\begin{equation}
    u_0 \in L^2(\mathbb{R}^N), \quad \quad \emptyset \ne \text{supp}(u_0)\subseteq [-R_0,R_0]^N
\end{equation} \noindent 
Then there is a solution $\tilde{u} \ne 0$ such that 
\begin{equation} \label{support}
\text{supp} (\tilde{u}(\cdot,t)) \subseteq \prod_{i=1}^N [-R_j(t),R_j(t)], 
\end{equation} \noindent for any $t<T$, where
\begin{equation} \label{Rj}
R_j(t)= 2R_0+Ct^{\frac{N(\bar{p}-p_j)+\bar{p}}{\lambda p_j}} ||u_0||_1^{\frac{\bar{p}}{p_j} \frac{p_j-2}{\lambda}}, \quad \quad \quad \lambda=N(\bar{p}-2)+\bar{p}.
\end{equation}
\end{theorem} \noindent
Moreover, they proved the following $L^\infty$-$L^1$ estimates of the decay for the solution.
\begin{theorem} Let $\bar{p}<N$ and let $u \in \cap_{i=1}^N \L^{p_i}(\Sigma_T)$ solve \eqref{anisocauchy} for $u_0 \in L^1(\R^N) \cap L^2(\R^N)$. Then if $p_i>2, \, \forall i=1,..,N$
the following estimate holds true for any $\tau \in [0,T]$
\begin{equation} \label{decay}
|| u(\cdot,t)||_{L^\infty(\mathbb{R}^N)} \leq C t^{-\frac{N}{\lambda}} ||u_0||_{L^1(\mathbb{R}^N)}^{\frac{\bar{p}}{\lambda}}.
\end{equation}
\end{theorem}

\subsection*{The anisotropic Fokker-Planck equation.}
We consider a similar continuous transformation as \eqref{transf-data}, owing the choiche of the right exponent to the decay of a solution to \eqref{anisocauchy}, and we perform the following formal calculations.
\begin{equation} \label{anisotransform}
    u(x,t)= t^{-\beta}v \bigg(x_1 t^{\alpha_1},..., x_N t^{\alpha_N},t \bigg)=t^{-\beta}v(y_{1},..,y_N,t), \quad 
    \begin{cases} y_i= x_i t^{\alpha_i} ,\\ \frac{\partial}{\partial x_i}=  t^{\alpha_i} \frac{\partial }{\partial y_i}.
    \end{cases}
\end{equation} \noindent
We calculate formally
\[
u_t=-\beta t^{-\beta-1}v+ t^{-\beta}\bigg[ \sum_{i=1}^{N} \bigg( \frac{\partial}{\partial y_i} v \bigg) \frac{\partial y_i}{\partial t}+v_t \bigg] =-\beta t^{-\beta-1}v+  t^{-\beta} \sum_{i=1}^{N} \bigg( \frac{\partial}{\partial y_i} v \bigg) \bigg[ \frac{\alpha_i x_i t^{\alpha_i}}{t}  \bigg] + t^{-\beta} v_t,\] being
\[
\frac{\partial}{\partial x_i} u= t^{\alpha_i-\beta} \frac{\partial}{\partial y_i} v.
\] \noindent
We substitute these into \eqref{prototype} to get
\[ -\beta t^{-\beta-1}v+  t^{-\beta} \sum_{i=1}^{N} \frac{\alpha_i y_i}{t} \bigg( \frac{\partial}{\partial y_i} v \bigg)     + t^{-\beta} v_t=\sum_{i=1}^{N} t^{\alpha_i} \frac{\partial}{\partial y_i} \bigg(t^{(\alpha_i-\beta)(p_i -1)} \bigg| \frac{\partial}{\partial y_i}v   \bigg|^{p_i-2}\frac{\partial}{\partial y_i}v  \bigg). \]
\noindent 
Re-ordering and multiplying each term for $t^{\beta+1}$ we get
\[
t v_t= \beta v - \sum_{i=1}^{N} \alpha_i y_i \frac{\partial}{\partial y_i} v +\sum_{i=1}^{N} t^{(\alpha_i-\beta)(p_i-1)+\alpha_i+ \beta+1} \frac{\partial}{\partial y_i} \bigg( \bigg| \frac{\partial}{\partial y_i}v   \bigg|^{p_i-2}\frac{\partial}{\partial y_i}v  \bigg)=
\]
\[
\beta v + \sum_{i=1}^{N} \alpha_i v+ \sum_{i=1}^{N}\frac{\partial}{\partial y_i} \bigg[ \bigg(\bigg| \frac{\partial}{\partial y_i}v   \bigg|^{p_i-2}\frac{\partial}{\partial y_i}v \bigg) -\alpha_i y_i v     \bigg],
\]
by choosing 
\[(\alpha_i-\beta)(p_i-1)+\alpha_i+ \beta+1=0,
\]
which means
\begin{equation} \label{alphai}
      \alpha_i= \beta- \frac{1+2\beta}{p_i} < 0.
\end{equation} 
This is an Euler equation. So, by redefining $v(y,t)= w(y,\ln (t))$ the equation \eqref{prototype} becomes the non-homogeneous Fokker-Planck equation
\begin{equation} \label{Eq:1:2}
     w_t= \bigg( \beta+\sum_{i=1}^{N} \alpha_i \bigg) w+ \sum_{i=1}^{N}\frac{\partial}{\partial y_i} \bigg[ \bigg(\bigg| \frac{\partial}{\partial y_i}w   \bigg|^{p_i-2}\frac{\partial}{\partial y_i}w \bigg) -\alpha_i y_i w     \bigg].
\end{equation} \noindent
If, according to \eqref{decay}, we consider
\begin{equation} \label{beta}
\beta= \frac{N}{N(\bar{p}-2)+\bar{p}},
\end{equation} then the equation reduces to
\begin{equation} \label{Eq:1:2bis}
     w_t= \sum_{i=1}^{N}\frac{\partial}{\partial y_i} \bigg[ \bigg(\bigg| \frac{\partial}{\partial y_i}w   \bigg|^{p_i-2}\frac{\partial}{\partial y_i}w \bigg) -\alpha_i y_i w     \bigg].
\end{equation} \noindent 
\begin{remark}
Equation \eqref{Eq:1:2bis} conserves the $L^1(\Omega)$-norm in time.\newline
Moreover, a solution to the stationary version of \eqref{Eq:1:2bis} would give us the wanted Barenblatt Fundamental solution to \eqref{prototype}.
\end{remark} \noindent
This anisotropic Fokker-Planck type equation is deeply different from its isotropic counterpart \eqref{Fk-Pk-isotropa}. Anisotropy does not permit the identification of a single variable ODE as in \eqref{monoFokker}, and this is physically evident and due to the lack of radial symmetry of the diffusion process in consideration: there is no homogeneous flux here to be vanished. Moreover the steady equation
\begin{equation} \label{Fk-PK-anisotropa-stazionaria}
\sum_{i=1}^{N}\frac{\partial}{\partial y_i} \bigg[ \bigg(\bigg| \frac{\partial}{\partial y_i}w   \bigg|^{p_i-2}\frac{\partial}{\partial y_i}w \bigg) -\alpha_i y_i w     \bigg], \quad \text{in} \quad \Omega \subset \R^N,
\end{equation} \noindent is not a variational one i.e. it is not known if it can be written as the Euler Lagrange equation of an energy functional. Moreover, its monotonicity and coercivity properties suffer heavily the second term influence relatively to the length in the $i$-th direction of the medium $\Omega$. These considerations leading to the difficulty of an explicit formula as in the previous case \eqref{Barenblatt}, the existence and the main properties characterizing a Barenblatt Fundamental solution may be derived by the simpler original equation \eqref{prototype} and then defining a suitable function which solves the steady Fokker-Planck equation \eqref{Fk-PK-anisotropa-stazionaria}. This would ensure that the solution to \eqref{prototype} found has the properties of Theorem \ref{Sunra}, which characterize a Barenblatt Fundamental Solution.

\subsection*{On the solvability of the anisotropic Cauchy Problem with measure initial data.} \label{CPTedeev}
We consider the prototype problem with measure initial data, i.e
\begin{equation} \label{CPanisotropic+MeasureData}
    \begin{cases}
    u_t - \sum_{i=1}^N (  |u_{x_{i}}|^{p_i-2} u_{x_i} )_{x_i} =0 , \quad (x,t) \in \R^N \times [0,T],\\
    u(x,0)=u_0(x), \quad x \in \R^N.
    \end{cases}
\end{equation} \noindent We begin the study of a weak solution to \eqref{CPanisotropic+MeasureData} i.e. a function $u \in C(0,T; L^1(\R^N)) \cap L^{\bf{p}}(0,T;W^{1,{\bf{p}}}(\R^N))$ such that for each open bounded $\Omega \subset \R^N$ and for all $t \in [0,T)$ satisfies for all test function $\varphi(x,t) \in W^{1, \infty}([0,T,L^{\infty}(\Omega)) \cap L^{\infty}([0,T],W^{1, \infty}_o(\Omega))$ the equality
\begin{equation} \label{weak-solution-2}
    \int_{\Omega} u \varphi(x,t) dx + \sum_{i=1}^N \int_0^t \int_{\Omega} |u_{x_i}|^{p_i-2} u_{x_i} \varphi_{x_i}  dx d\tau= \int_{\Omega} \varphi(x,0) d u_{0}+ \int_0^t \int_{\Omega} u \varphi_{\tau}(x, \tau) dx d\tau.
\end{equation} \noindent
This has been done in \cite{Tedeev-CP}, \cite{Tedeev-estimates} for more general doubly nonlinear anisotropic equations. We recall the notation $\lambda= N(\overline{p}-2)+\overline{p}$. In \cite{Tedeev-estimates} the authors prove a generalised version of the following a priori estimates.
\begin{theorem} \label{a-priori-estimates}
Consider the problem \eqref{CPanisotropic+MeasureData} with $2<p_i \leq {\bar{p}}\bigg(1+\frac{1}{N} \bigg)$, $u_0(x) \ge 0$ and 
\begin{equation} \label{MAGIC-CONDITION}   
\| | u_0 \||_r:= \sup_{\rho \ge r} \rho^{-\frac{\lambda}{N}}  \int_{B_{\rho}} u_0(x) dx < \infty,\quad  r>0,
\end{equation} \noindent being
\[  
\quad B_{\rho}:= \bigg{\{} x \in \R^N| |x_i| \leq \frac{\rho^{\frac{{\bar{p}}(p_i-2)}{p_i({\bar{p}}-2)}}}{2} \bigg{\}}.
\]
Define by monotonicity $M_{\infty}:= \lim_{r \rightarrow \infty} |\| u_0 \||_r$ and for a $\gamma>0$ to be specified later \begin{equation} \label{waiting-time}
T_{*}:= \begin{cases} 
\infty, & \text{if} \quad M_{\infty}=0,\\
\bigg( \frac{M_{\infty}}{\gamma} \bigg)^{\frac{N(\overline{p}-p_N)+\overline{p}}{\overline{p}(p_N-2)}}, & \text{if} \quad M_{\infty}\ge \gamma,\\
\bigg( \frac{M_{\infty}}{\gamma} \bigg){\frac{N(\overline{p}-p_1)+\overline{p}}{\overline{p}(p_1-2)}}, & \text{if} \quad M_{\infty}< \gamma.
\end{cases}
\end{equation} \noindent 
Then there exists a positive constant $\gamma(p_i,N)$ such that every nonnegative weak solution to \eqref{CPanisotropic+MeasureData} defined on $[0,T_{*}]$ must satisfy the following estimates for all $ t, {\bar{t}} \in (0,T_{*})$:
\begin{equation} \label{estimate1}
    |\| u(\cdot,t) \| |_r \leq C |\| u_0 \||_r,
\end{equation}
\begin{equation} \label{estimate2}
    \| u(\cdot,t) \|_{L^{\infty}(B_r)} \leq C r^{\frac{\overline{p}}{N}} t^{- \frac{N}{\lambda}} \||u_0\||_r^{\frac{\overline{p}}{\lambda}}, 
\end{equation}
\begin{equation} \label{estimate3}
    \sum_{i=1}^N \int_0^t \int_{B_r} |u_{x_i}|^{p_i-1} dx d\tau <C(r,t),
\end{equation}
\begin{equation} \label{estimate4}
    \sum_{i=1}^N \int_{\overline{t}}^t \int_{B_r} |u_{x_i}|^{p_i} dx d\tau <C(r,t, \overline{t}).
\end{equation}
\end{theorem} \noindent 
\begin{remark} For $p_i=p,\, \forall i=1,..,N$ estimates \eqref{estimate1}, \eqref{estimate2}, \eqref{estimate3}, \eqref{estimate4} and the number $T_{*}>0$ do coincide with the ones of Section \ref{SolvingCP} for the isotropic equation found in \cite{DB-He}. Secondly, it is interesting to observe that the lifetime of the solution is determined by the largest exponent $p_N$ in case of large initial mass $\|u_0\||_r$ while it is determined by the smaller $p_1$ in case of a modest initial mass.
\end{remark} \noindent 
\section*{Future strategy and conclusion.}
In this note we have proven the strong connection between the Barenblatt Fundamental solution and the solutions to the stationary equation \eqref{Fk-PK-anisotropa-stazionaria}. We have shown the existence of solutions to \eqref{Eq:1:2bis} thanks to a recent result in \cite{Tedeev-CP}. However, this is not enough to use this result to prove regularity results. Indeed, we can invoke the previous Theorem to find a solution $u$ to \eqref{prototype}. We already know that there exists a solution of $u$ that satisfies the growths \eqref{support}, \eqref{decay}. But what is missing, to repeat the same ideas of Section \ref{harnackpotente}, is a nice description {\it from below} of the support of $u$. The aim of our next papers is to carry on a deep analysis of the interplay between these two equations and to develop the necessary tools for deriving regularity results and Harnack inequalities for nonnegative solutions to \eqref{prototype}.


\end{document}